\input amssym.def
\input  amssym.tex

\documentstyle[11pt]{article}

\title{Natural Internal Forcing Schemata extending ZFC. Truth in the Universe?}

\author{Garvin Melles\thanks{Would like to thank Ehud Hrushovski
for supporting him with funds from NSF Grant DMS 8959511}\\Abraham
Fraenkel Center for Mathematical Logic 
\\Bar Ilan University\\Institute of Mathematics\\Hebrew University of Jerusalem}

\newcommand{\proof}{{\sc proof} \hspace{0.1in}}

\newcommand{\sub}{\subseteq}
\newcommand{\al}{\alpha}
\newcommand{\be}{\beta}
\newcommand{\ga}{\gamma}

\begin{document}
\mathsurround=.1cm
\maketitle

\begin{center}
{\bf INTRODUCTION}
\end{center}

Mathematicians are one over on the physicists in that they already have a
unified theory of mathematics, namely set theory. Unfortunately the
plethora of independence results since the invention of forcing has
taken away some of the luster of set theory in the eyes of many
mathematicians. Will man's knowledge of mathematical truth
be forever limited to those theorems derivable from the standard
axioms of set theory, $ZFC?$ This author does not think so, he feels
that set theorists intuition about the universe is stronger than ZFC.
Here in this paper, using part of this intuition, he introduces some
axiom schemata which he feels are very natural 
candidates for being considered as part of the axioms of set theory. 
These schemata assert the existence of many generics over simple
inner models. The main purpose of this article is to present 
arguments for 
why the assertion of the existence of such generics belongs to the
axioms of set theory. 

Our central
guiding principle in justifying the axioms is what Maddy called the
rule of thumb maximize in her survey article on the axioms of set
theory, [BAI] and [BAII]. More specifically, our intuition conforms
with that expressed by Mathias in his article ``What is Maclane
Missing?'' challenging Mac Lane's view of set theory. 
\begin{quotation}
This might be a good moment to challenge one of Mac Lane's opinions,
which I believe to rest on a misconception. On page 359 of his book he
writes, after reflecting on the plethora of independence results, that
``for these reasons 'set' turns out to have many meanings, so that the
purported foundations of all Mathematics upon set theory totters.''
Elsewhere, on page 385, he remarks that ``the Platonic notion that there
is somewhere the ideal realm of sets, not yet fully described, is a
glorious illusion.''

I would suggest a contrary view: independence results within set
theory are generally achieved either by examining an inner model of
the universe (an inner model being a transitive class containing all
ordinals) or by utilizing forcing to obtain a larger universe of which
the original one is an inner model. The conception that begins to seem
more and more reasonable with the advance of the inner model program
on the one hand and a deeper understanding of iterated forcing on the
other is that within one enormous universe there are many inner
models, and the various ``independence arguments'' may be reworked to
give positive information about the way various inner models  relate
to one another. Far from undermining the set theoretic point of view,
the various techniques available for building models actually promote
that unity.
\end{quotation}

One of this author's reasons for having an
intuition about sets similar to that expressed by Mathias is given by a very short look back at 
the history of the development of mathematics. Mathematics began
with the study of mathematical objects very physical and concrete in
nature and has progressed to the study of
things completely imaginary and abstract. Most mathematicians now accept these objects
as mathematically legitimate as any of their more concrete
counterparts. It is enough that these objects are consistently
imaginable, i.e., exist in the world of
set theory.  Applying the same intuition to set theory itself, we get
some motivation for why we
should accept as sets generic objects over inner models. 

Using the rule
of thumb maximize, the principle of reflection and esthetics as a basis we will
now procede to give further more concrete arguments which 
say that theories $T$ extending $ZFC$ implying the existence of many
generics over inner models reflect truth about the universe of sets
$V.$ While the arguments are not technically sophisticated, the author
hopes this will not detract from the axioms intuitive appeal.

\pagebreak 

\begin{center}
{\bf SET THEORETIC PRELIMINARIES}
\end{center}

In order to make this article accessible to the
general reader, we review in this section some of the set theoretic basics needed to
understand the arguments in favor of the schemata. As this review is
sparse, any gaps in the readers understanding can be filled in with
reading the appropriate sections of [Jech2] or [Kunen]. Set theorists can
safely skip this section except for the last three definitions, which 
are not standard.

A formula in the language of set theory is a formula in the first
order predicate claculus built from atomic formulas of the form $x=y$
and $x\in y.$ We denote the universe of sets by $V$ and the
satisfaction relation by $\models.$ A class is a
collection of sets satisfying some formula of set theory. (So all sets
are classes.)  We denote the class of ordinals by $Ord.$  A set
$Y\subseteq X$ is a definable subset of $X$ if for some formula of set
theory $\varphi(x),$
$$\forall y(y\in Y\ \leftrightarrow\ X\models\varphi(y)\,)$$
If $X$ is a set, then
$Def(X)$ is the set of all definable subsets of $X$ and $P(X)$ the set
of all subsets of $X.$ 
If we assume the axiom of foundation then the universe of sets $V$ can be
written as
$$V=\bigcup_{\alpha\in Ord}V_{\alpha}$$
where $V_{\alpha+1}=V_{\alpha}\cup P(V_{\alpha})$ and if $\alpha$ is a
limit ordinal then 
$$V_{\alpha}=\bigcup\limits_{\beta<\alpha}V_{\beta}$$
By a real we will usually mean a subset of $\omega,$ but within set theory a real
can also mean an element of $V_{\omega+1},\ \omega^{\omega}$ or 
${\Bbb R}.$  The constructible universe $L$ (the class of
constructible sets) can be written as 
$$L=\bigcup_{\al\in Ord}L_{\al}$$
where $L_{\emptyset}=\emptyset,$ $L_{\al+1}= L_{\al}\cup Def(L_{\al})$ and if ${\al}$ is a limit
ordinal then $L_{\al}=\bigcup\limits_{\be < \al}L_{\be}.$ $L$ is a
model of $ZFC\ +\ GCH.$ $L$ is absolute in the sense that if $M$ and
$N$ are transitive models of $ZF$ with the same ordinals, then
$L^M=L^N,$ i.e., the class of things $M$ thinks is constructible is the
class of things which $N$ thinks is constructible. Similarly, if $X$
is a transitive set,
$L(X)$ the class
of sets constructible from $X,$ can be written as 
$$L(X)= \bigcup_{\al\in Ord}L(X)_{\al}$$ 
where $L(X)_{\emptyset}=X,$ $L(X)_{\al+1}=
L(X)_{\al}\cup Def(L(X)_{\al})$ and if $\al$ is a limit ordinal, then
$L(X)_{\al}=\bigcup\limits_{\be <\al}L(X)_{\be}.$ In general $L(X)$ is
a model of $ZF$ but not of $AC.$ If $M$ is a transitive model and
$x\in M,$ we say  $x$ is definable in $M$ if for some formula
$\varphi(v)$ of set theory, $M\models \exists !v\varphi(v)\ \wedge\
\varphi(x).$ 
    Now we turn to a short review of forcing. If $P$ is a partial
order, a subset $D$ of $P$ is said to be dense if for every $p\in P,$
there is a $d\in D$ such that $d\leq p.$ If $p,q\in P$ such that there
is no $r\in P$ such that $r\leq p$ and $r\leq q,$ then we say $p$ and $q$ are incompatible,
written, $p\perp q.$ If $P$ is a partial
order in a transitive model $M$ of $ZFC,$ then a subset $G$ of $M$ is
said to be $M$ generic if $G$ intersects every dense subset of $P$
in $M.$  If $M$ is a transitive
model of $ZFC$ and $P$ is a partial order in $M$ such that for every $p\in
P$ there is $r$ and $q$ in $P$ such that $r\leq p,\ q\leq p,\ r\perp
q,$ and if $G$ is an
$M$ generic subset of $P,$ then there is a transitive model $M[G]$ of $ZFC$
such that 
\begin{enumerate}
\item $G\in M[G]$ and $G\not\in M$  
\item $Ord^M=Ord^{M[G]}$
\item If $N$ is a transitive model such that $M\sub N,$ 
$Ord^M=Ord^N,$ and $G\in N,$ then $M[G]\sub N.$ 
\end{enumerate}
$\!M[G]$ is called a forcing extension of $M.$ A partial order $P$ is
separative if for every $p,q\in P, $ $p\not\leq q\ \rightarrow\
\exists r(r\leq p\ \wedge\ r\perp q).$ Without loss of
generality we can assume all the partial orders we use are separative. 
Associated with every separative partial order is the
set of regular open subsets of $P,$ denoted $r.o.(P).$

\vspace{.15in}

\noindent $r.o.(P)=
\Big\{S\subset P\mid \forall p,q\in P\ p\in S\ \wedge\ q\leq p\
\rightarrow\ q\in S\ \wedge $
\begin{flushright}
$\forall p\in P\big(\forall r\in P(r\leq p \ \rightarrow\ \exists q\in S(q\leq r)\,)\ 
\rightarrow\ p\in S\big)\Big\}$
\end{flushright}

\vspace{.05in}

\noindent With the appropriate interpretation of $+$ and $\bullet$, $r.o.(P)$
is a complete Boolean algebra. If $\kappa$ and $\lambda$ are
cardinals, a Boolean algebra $B$ is $(\kappa,\lambda)$
distributive iff every collection of $\kappa$ partitions of $B$ of size at
most $\lambda$ has a common refinement.

\vspace{.1in}

\noindent{\bf Theorem 1.}
If $M$ is a transitive model of $ZFC,$ $\kappa$ is a cardinal in $M$
and $M[G]$ is a forcing extension of $M$ via the partial order $P,$
then $M[G]$ has no functions
$f:\kappa\rightarrow\kappa$ not in the ground 
model if and only if $r.o.(P)$ is $(\kappa,\kappa)$ distributive.\\
\proof See [Jech2].

\vspace{.1in}

\noindent{\bf Definition 1.}
A subset $r$ of $\omega$ is said to be absolutely definable if for
some $\Pi_1$ formula $\theta(x),$

\begin{enumerate}
\item $V\models \theta(r)$
\item $V\vdash \exists!\,x\theta(x)$

\end{enumerate}  
A canonical example of an absolutely definable real is $0^{\#}.$ 

\vspace{.1in}

\noindent{\bf Definition 2.}
$x\in V$ is said to be weakly absolutely definable if for some
$\Sigma_1$ formula $\psi(x)$ 
$$V\models \forall y(y\in x\ \leftrightarrow\ \psi(y)\,)$$

\noindent{\bf Definition 3.}
$x\in V$ is said to be weakly absolutely definable of
the form $V_{\alpha}$ if for some ordinal $\alpha$ definable in $L,$ 
$$V\models \forall y(y\in x\ \leftrightarrow\ \rho(y)\leq \alpha)$$ 

\vspace{.1in}

\begin{center}
{\bf THE SCHEMATA}
\end{center}
 
In this section we list the Schemata. Schemata we know how to prove the
consistency of assuming the existence of a countable transitive model
of $ZFC$ we label with a (*). The consistency of the other schemata
are just conjectured, and we
write (Conj) besides those. Schemata which follow from large cardinal
assumptions we write (FLC) besides. Note that the definitions are informal as
the formal versions are unwieldy.

\vspace{.1in}
 
\noindent{\bf Definition 4.} (*) (FLC) $IFS(L)$ is the axiom schema which says
for every formula $\phi(x),$ if
$L\models$ there is a unique partial order $P$ such that $\phi(P),$
then there is a $L$ generic subset of $P$ in the universe $V.$

\vspace{.1in}

\noindent{\bf Definition 5.} (*) (FLC) $IFS(L[r])$ is the axiom schema of set
theory which says if $r$
is an absolutely definable real then every partial order $P$ definable in
$L[r]$ has an $L[r]$ generic subset.

\vspace{.1in}

\noindent{\bf Definition 6.} (*) $IFS(L({\Bbb R})$ is the axiom schema of
set theory which says if $P$ is a partial order definable in
$L({\Bbb R})$ such that 
$$\Vdash {\Bbb R}={\Bbb R}^{V^P}$$
(By $\Vdash$ we mean $\Vdash$ in V) then there exists a $L({\Bbb R})$
generic subset $G$ of $P.$ 

\vspace{.1in}

\noindent{\bf Definition 7.} (*) $IFS(L(V_{\alpha}))$ is the axiom schema
which says if $V_{\alpha}$ is weakly absolutely definable of the form
$V_{\alpha}$ and $P$ is a partially ordered set definable in
$L(V_{\alpha})$ such that 
$$\Vdash V_{\alpha}=V_{\alpha}^{V^P}$$
then there is a $L(V_{\alpha})$ generic subset $G$ of $P.$

\vspace{.1in}

\noindent{\bf Definition 8.} (Conj) $IFS(\forall L(V_{\alpha}))$ is the
axiom schema which says for all $\alpha\in Ord,$ if $P$ is a partially orderd
set definable in $L(V_{\alpha})$ such that
$$\Vdash V_{\alpha}=V_{\alpha}^{V^P}$$
then there is a $L(V_{\alpha})$ generic subset $G$ of $P.$

\vspace{.1in}

\noindent{\bf Definition 9.} (*) $IFS^-(\forall L(V_{\alpha}))$ is the
axiom schema which says for each ordinal $\alpha$ if $P$ is a partial
order definable in 
$L(V_{\omega+\alpha})$ such that $P$ is $\aleph_{\beta}$ closed for
each $\beta<\alpha,$ then there is a $L(V_{\omega+\alpha})$ generic
subset $G$ of $P.$

\vspace{.1in}

\noindent {\bf Definition 10.} (Conj) $IFS$ is the axiom schema of set theory
which says for every
weakly absolutely definable set $X,$ for every partial order 
$P$ definable in $L(X),$ if
$$\Vdash\ X^{V^P}=X$$
then there exists an $L(X)$ generic subset $G$ of $P.$

\vspace{.1in} 

\noindent {\bf Definition 11.} (FLC\ +\ CH) $IFS\restriction X\subseteq {\Bbb R}$
makes the same claim as $IFS$ but only for $X\subseteq {\Bbb R}.$

\vspace{.1in}

\noindent {\bf Definition 12.} (*) $CIFS$ is the axiom schema of set
theory which says that for each regular cardinal $\aleph_{\alpha}$ all
definable subsets of $L(V_{\omega+\alpha})$ are of size at most
$\aleph_{\al}.$ 

\vspace{.1in} 

\noindent We conjecture that $IFS(\forall L(V_{\alpha}))$ is
consistent, since it is a natural generalization of
$IFS(L(V_{\alpha})).$ 
The intuition behind the conjecture that $IFS$ is consistent is
somewhat more nebulous. One can look upon $IFS(L(V_{\alpha}))$ as
saying the universe has a sort of minimal largeness with respect to
the $L(V_{\alpha})$ because it implies there are many $L(V_{\al})$ generics in the
universe. The $L(V_{\alpha})$ provide a reference frame
from which to measure the size of the universe, since $L(V_{\alpha})$
is absolute for any class model of $ZFC$ containing $V_{\alpha}$ and
$L(V_{\al})$ generic subsets for partial orders definable in $L(V_{\al})$
maintain their genericity under extensions as long as $V_{\al}$ (and
the class of ordinals) is not
changed. If
$X$ is weakly absolutely definable, $L(X)$ is also absolute for any
class model of $ZFC$ containing $X.$ So $IFS$ is a natural
generalization of $IFS(L(V_{\alpha}))$ 
implying the universe has a minimal largeness with respect to each of the $L(X).$
For a given weakly absolutely definable set $X$ the consistency of
$IFS\restriction X$ is
easy to show.

\vspace{.1in}

\vspace{.1in}

\begin{center}
{\bf WHY SHOULD THE SCHEMATA SHOULD HOLD IN $V$?}
\end{center}

Our version of the 'rule of thumb maximize' will take the form of the
following three principles:
\begin{enumerate}
\item $V\models ZFC$
\item $V$ is large with respect to $Ord$
\item $V$ is large with respect to each of the $V_{\al}$
\end{enumerate}
In order to get a better handle on what principles two and three mean,
we shall use countable transitive sets  as models for transitive
classes scaled down to a countable size. We will take a look at
countable transitive models of $ZFC$ satisfying principles two and
three and look for common and esthetically pleasing properties among
them, i.e., properties that we think $V$ itself should satisfy. 
In order to see what principle two gives us, we fix the height of the
models under consideration i.e., we assume all our models have the
same set of ordinals $\alpha.$ (And we also assume of course that
countable transitive models of $ZFC$ with ordinals $\alpha$ exist.)   
So we are using $\al$ as a model for $Ord.$ Now there is a unique
countable transitive model of $ZFC+V=L$ with 
height $\alpha,$ namely the set $L_{\alpha}.$ $L_{\alpha}\subseteq M$
for every $M$ which is a countable transitive model of $ZFC$ of height
$\alpha$ and $L^M=L_{\al}.$ So the statement $V=L$ expresses  a kind of minimal
property, the opposite of what we are looking for. On the other hand, the statement
$IFS(L)$ is a kind of minimal maximality condition among the countable
transitive models of $ZFC$ with given height. Why? Suppose $M,N$ are
countable transitive models of $ZFC$ such that $M\models
IFS(L)$ with $M\sub N$ and $Ord^N=Ord^M.$ Then $N\models IFS(L)$
since the interpretation of $L$ and of the $L$ generics for the
various partial orders definable in $L$ are absolute. Furthermore,
larger models tend toward $IFS(L),$ i.e., given any countable transitive
model $M$ of $ZFC$ and any finite list $P_1,\ldots,P_n$ definable in
$L^M,$ if we let $N$ be the forcing extension of $M$ by the partial
order $P_1\times\ldots\times P_n$ then $N$ has the same height as $M$
and satisfies $IFS(L\restriction\{P_1,\ldots,P_n\}).$ So $IFS(L)$ is a
natural closure condition on the countable transitive models of $ZFC$ of given
height. The arguments for the axiom schema $IFS(L[r])$ have similar justifications. 
As we consider that the relationships among countable transitive
models of $ZFC$ are reflections of the relationships among transitive
class models of $ZFC,$ we argue that $IFS(L)$ and $IFS(L[r])$ should
hold in $V.$ If $V\not\models IFS(L),$ it would be as if the universe
had an artificial boundary. It seems it would be an artificial constraint on
$V$ if for some $P$ a partial order definable in $L$ there is no
$L$ generic. Note that under $ZFC\ +\ IFS(L[r]),$ all
the generics asserted to exist by the axioms of $IFS(L[r])$ are in
$L({\Bbb R})=L(V_{\omega+1}),$ so $IFS(L[r])$ is really a
schema about the structure of $L({\Bbb R}).$ Note also that to be more
formal and to work strictly within $ZFC$ we could have made our
arguments using countable transitive models of arbitrarily large
finite parts of $ZFC$ and the schemata. 

Why do we work with partial orders $P$ definable in $L$ and not all
$P\in L?$ In the first place axioms asserting the existence of
generics for all $P\in L$ are inconsistent with $ZFC,$ but the main
point is that we are interested not in countable transitive models but
in proper class models of $ZFC$ and forcing only gives the relative
consistency of extensions of $L$ of the form $L[G]$ only for those
$G\subseteq P$ where $P$ is a partial order definable (without
parameters) in $L.$ In
keeping with our principal of maximality we reinterpret this to mean
that such generic extensions of $L$ actually exist.  

To investigate the consequences of principle three, we fix both the
height and the width at stage 
$\omega+1$ (i.e., $V_{\omega+1}$) among the
models (which we can assume satisfy $ZFC+IFS(L[r])$) under
consideration. (All transitive models of $ZFC$ have the 
same first $\omega$ stages in the cumulative heirarchy.) 
So we are using some countable ordinal as a model for $Ord$ and some
countable set of reals as a model for the reals. Arguing as before we see that
$IFS(L({\Bbb R}))$ is a natural closure condition on this class of models,
implying a minimal kind of maximality. Similarly, we argue that among the
countable transitive models of given height and set of reals
satisfying $ZFC+IFS(L[r])$ a natural closure property is that all sets
definable in $L({\Bbb R})$ are of size at most $\aleph_1$ since the
canonical forcing which collapse definable elements of $L({\Bbb R})$
to size $\aleph_1$ are $\omega$ closed and therefore do not add reals.
Continuing to make use of our third principle, similar reasoning works for all 
the definable stages $V_{\alpha}$ so we are lead to
$IFS(L(V_{\alpha})),\ IFS^-(\forall L(V_{\al})),\ IFS(\forall
L(V_{\al})),$ and $CIFS.$ 

Another justification for the schemata (see [BA I] page 492-493) is that they are
a way of making the power set thick. More precisely, insteads of making
$P(V_{\alpha})-V_{\alpha}$ large, they make
$L(V_{\alpha+1})-L(V_{\alpha})$ large, a slight variant of
the notion that the power set operation should be large. This is one
of the appeals behind $IFS(\forall L(V_{\al})),\ IFS^-(\forall
L(V_{\al}))$ and $CIFS.$ Even under
$IFS^-(\forall L(V_{\alpha})),$ for all
regular cardinals $\aleph_{\alpha},$ 
$L(V_{\omega+\alpha+1})-L(V_{\omega+\alpha})\neq \emptyset.$

\vspace{.1in}

\begin{center}
{\bf CONNECTIONS WITH LARGE CARDINALS}
\end{center}

It is not hard to see that $0^{\#}$ exists implies $IFS(L)$ since as
we show later that $IFS(L)$ is equivalent to the assumption that every
set definable in $L$ is countable. So $IFS(L)$ is a kind of intrinsic
support for the large cardinal axiom $0^{\#}$ exists.  
Similarily the picture of the universe given by $IFS(L[r])$ is related to that
under the assumption of a measurable cardinal. If a measurable
cardinal exists than $r^{\#}$ exists for every $r\subseteq \omega,$ so
that means for every $r\subseteq \omega,$ every set definable in
$L[r]$ is countable. As we shall soon prove, $IFS(L[r])$ holds if and
only if
for every $r$ which is an absolutely definable real, every set definable in $L[r]$ 
is countable. So again $IFS(L[r])$ provides a kind of intrinsic support for
large cardinal axioms, in that they give at some level similar pictures of the universe,
even though the consistency strength of the large cardinal axioms are
much greater than that of $IFS(L)$ or $IFS(L[r]).$ The most important
connection between large cardinals and the schemata known to the author is the fact which
was pointed out to him by Woodin that under large cardinal hypotheses, 
$IFS(\restriction X\subseteq {\Bbb R})$ is
equivalent to $CH.$

\vspace{.1in}

\begin{center}
{\bf CONSISTENCY FROM A COUNTABLE TRANSITIVE MODEL OF ZFC}
\end{center}

\noindent {\bf Theorem 2.}
Let $\langle\theta_i\mid i<n\rangle$ 
and $\big\{\varphi_{ij}(x)\mid i<n, j<m,\big\}$ be  finite sets  of
formulas with the $\theta_i$ being $\Pi_1.$  Let $M$ be a countable
transitive model of $ZFC.$ Then there exists a countable transitive
model $M'$ of $M$ with the same ordinals as $M$ such that for each
$i<n,$ 
$$M'\models \exists !r\theta_i(r)\ \wedge\ \theta_i(r_i)\
\rightarrow$$
$$\bigwedge_{j<m}\Big(\,(\,L[r_i]\models\exists !P(\varphi_{ij}(P))\
\wedge\ \varphi_{ij}(P_i)\,)\ \rightarrow\ \exists G\sub P_{ij}(G\hbox{
is }L[r_i]\hbox{ generic})\Big)$$
\proof Let $\alpha^*\in M$ such that 
$$M\models \alpha^*> sup\Big\{|{\cal D}_{ij}|\mid \exists
!x\theta_i(x)\ \wedge\ \theta(r_i)\ \wedge\ r_i\sub \omega\ \wedge$$
$$L[r_i]\models \exists !P(\,\varphi_{ij}(P)\ \wedge\
\varphi_{ij}(P_{ij})\ \wedge\ {\cal D}_{ij}\hbox{ is the set of
dense subsets of }P_{ij})\Big\}$$
Let $P$ be the set of finite partial one to one functions from
$\alpha^*$ to $\omega.$ Let $M'=M[G]$ where $G$ is a $M$ generic
subset of $P.$ Note that by the Levy-Shoenfield absoluteness lemma,
if $M\models \theta_i(r_i)$ then also $M'\models \theta_i(r_i).$ Since
all the ${\cal D}_{ij}$ are countable in $M'$ the $P_{ij}$ have
$L[r_i]$ generic subsets in $M'.$
To finish the proof it is enough to prove the following claim.

\vspace{.1in}

\noindent Claim: If a formula $\psi(x)$ defines a real in $M[G]$ then it
is in $M.$\\
\proof Suppose $r$ is the unique real satisfying $\psi(x)$ in $M[G].$ 
Since $P$ is separative, if $p\in P$ and $\pi$ is an
automorphism of $P,$ then by [Jech 2] lemma 19.10, for every formula $\varphi(v_1,\ldots,v_n)$
and names $x_1,\ldots,x_n$ 
$$*\ \ \ \ p\Vdash \varphi(x_1,\ldots,x_n)\ \hbox{ iff }\ \pi p\Vdash
\varphi(\pi x_1,\ldots,\pi x_n)$$
Let $\varphi(x)=\exists\,Y(\psi(Y)\ \wedge\ x\in Y).$ Let $n\in
\omega.$ We will show that $||\varphi(\check{n})||=0$ or
$||\varphi(\check{n})||=1.$ If for no $p\in P$ does $p\Vdash
\big|\big|\varphi(\check n)\big|\big|$ then
$\big|\big|\varphi(\check n)\big|\big|=0.$ So let $p\in P$ such that $p\Vdash
\big|\big|\varphi(\check n)\big|\big|.$ By $*$ if $\pi$ is an automorphism
of $P$ then  $\pi p\Vdash
\big|\big|\varphi(\check n)\big|\big|.$ Let $\pi$ be a permutation of
$\omega.$ $\pi$ induces an automorphism of $P$ by letting for $p\in P,$
$dom\,\pi p=dom\,p$ and letting $\pi p(\alpha)=\pi(p(\alpha)).$ By
letting $\pi$ vary over the permutations of $\omega$ it follows that
$\big|\big|\varphi(\check n)\big|\big|=1.$ Let $\dot r$ be the name
with domain $\big\{\check n\mid n<\omega\big\}$ and such that
$$\dot r(\check n)=\big|\big|\varphi(\check n)\big|\big|$$ 
$i_G(\dot r)=r,$ but then
$r=\big\{n\mid \big|\big|\varphi(\check n)\big|\big|=1\big\}$ which
means it is in $M.$

\vspace{.1in}

\noindent {\bf Corollary 3.}
If there is a countable transitive model of $ZFC$ then 
$$Con(ZFC\ +\ IFS(L[r])\,)$$

\vspace{.1in}

\noindent {\bf Corollary 4.} $ZFC+IFS(L[r])\ +$'there are no
absolutely definable 
non-constructible reals' is
consistent. (Relative to the assumption of a countable transitive
model of $ZFC$)

\vspace{.1in}

\noindent {\bf Theorem 5.} If there is a countable transitive model of $ZFC$  then  
$$Con(ZFC+IFS(L(V_{\alpha}))\,)$$

\noindent \proof Let $\langle\theta_i\mid i<n\rangle$ 
and $\big\{\varphi_{ij}(x)\mid i<n, j<m,\big\}$ be  finite sets  of
formulas. Let $M$ be a countable transitive model of $ZFC.$ Without
loss of generality we can assume there exists ordinals
${\alpha_0,\ldots\alpha_{n-1}}$ such that 
$$L^M\models \exists!\alpha\theta_i(\alpha)\ \wedge\ \theta_i(\alpha_i)$$
and $\alpha_j<\alpha_k$ for $j<k.$ It is enough to find a forcing
extension $N$ of $M$ such that for each $i<n$ and $j<m$ for some
partial order $P_{ij}\in N$
$N\models$
$$L(V_{\alpha_i})\models \exists !x\psi_{ij}(x)\ \wedge\
L(V_{\alpha_i})\models \psi_{ij}(P_{ij})$$
$$\wedge \Vdash V_{\alpha_i}=V_{\alpha_i}^{P_{ij}}$$
$$\rightarrow\ \exists G(G\hbox{ is a }L(V_{\alpha_i})\hbox{ generic
subset of }P_{ij})$$
We define by induction on the lexicographical order of $n\times m$
sets $G_{ij}.$ Suppose $P_{ij}$ is a partial order 
definable in $L(V_{\alpha_i})^{M[\{G_{h,l}|h\leq
i,l<j\}]}$ by $\varphi_{ij}(x)$ and
there exists a $M[\{G_{h,l}|h\leq
i,l<j\}]$ generic subset of $P_{ij}$ not increasing 
$$V_{\alpha_i}^{M[\{G_{h,l}|h\leq
i,l<j\}]}$$ 
Then let $G_{ij}$ be such a $M[\{G_{h,l}|h\leq i,l<j\}]$
generic subset of $P_{ij}.$ (If not, let $G_{ij}=\emptyset.$) Let 
$$N=M[\{G_{ij}|i<n,j<m\}]$$

\vspace{.1in}

\noindent {\bf Theorem 6.} If there is a countable transitive model of $ZFC,$ then 
$$Con(ZFC+IFS(L(V_{\alpha}))+IFS(L[r])\,)$$\\
\proof Similar, just start with a model of enough of $IFS(L[r]).$

\vspace{.1in}

\noindent {\bf Theorem 7.} If there is a countable transitive model of
$ZFC$ then 
$$Con(ZFC+IFS^-(\forall L(V_{\alpha}))\,)\ \wedge\ Con(ZFC+CIFS)$$\\
\proof See the companion paper.

\vspace{.1in}

\begin{center}
{\bf SOME CONSEQUENCES AND SOME NICER FORMS}
\end{center}
Below we give some consequences and equivalents assuming $ZFC$ holds
in $V.$

\vspace{.1in}

\noindent {\bf Theorem 8.} $IFS(L({\Bbb R}))\vdash CH$\\
\proof Every bijection between a countable ordinal and a subset of
${\Bbb R}$ is an element of $L({\Bbb R})$ and
$\omega_1=\omega_1^{L({\Bbb R})}.$ So if $P=$ the set of bijections from countable ordinals into
${\Bbb R}$ then $P$ is a definable element of $L({\Bbb R}).$ Since $P$
is $\sigma$ closed, a $P$ generic over $V$ will not add any reals, so
by $IFS(L({\Bbb R}))$ there is a $G\subseteq P$ which is $L({\Bbb R})$ generic. 
If $\alpha$ is an ordinal
less than $\omega_1$ and $r$ is a real, let $D_{\alpha}=\big\{p\in
P\mid \alpha\in dom\,p\big\}$ and $D_r=\big\{p\in
P\mid r\in ran\,p\big\}.$ For each $\alpha<\omega_1,$ $G\cap
D_{\alpha}\neq\emptyset$ and for each $r\in
{\Bbb R},$ $G\cap D_r\neq\emptyset,$ so $\bigcup G$ is a bijection from
$\omega_1$ to ${\Bbb R}.$

\vspace{.1in}

\noindent {\bf Theorem 9.} $IFS(L({\Bbb R}))$ iff every $P$ definable
in $L({\Bbb R})$ such that $r.o.(P)$  is $(\omega,\omega)$ distributive
has an $L({\Bbb R})$ generic subset.\\
\proof By Theorem 1.

\vspace{.1in}

\noindent {\bf Theorem 10.} $IFS(\forall L(V_{\al}))$ iff every $\alpha\in
Ord$ and $P$ definable
in $L(V_{\omega+\al})$ such that $r.o.(P)$  is $(\aleph_{\be},\aleph_{\be})$ distributive
for each $\be<\al,$ has an $L(V_{\al})$ generic subset.\\
\proof By Theorem 1.

\vspace{.1in}

\noindent {\bf Theorem 11.} $IFS(L)\ \leftrightarrow\ $every set definable
in $L$ is countable. \\
\proof Certainly if every set definable in $L$ is countable, then
if $P$ is a partial order definable in $L$ then so is ${\cal D}$ the
set of dense subsets of $P$ in $L,$ so ${\cal D}$ is countable and therefore 
$P$ has a generic subset over $L$ in the universe. In the other
direction, if $s$ is a set definable in $L,$ then so is the partially
ordered set consisting of maps from distinct finite subsets of $s$ to distinct finite
subsets of $\omega,$ so a $L$ generic subset over the partial ordering
is a witness to $\big|s\big|=\omega.$

\vspace{.1in}

\noindent So $CIFS$ is a natural generalization of $IFS(L).$ 

\vspace{.1in}

\noindent {\bf Theorem 12.} $IFS(L[r])\ \leftrightarrow\ $for every
absolutely definable 
real $r,$ every set definable in $L[r]$ is countable.\\
\proof Similar to the previous proof.

\vspace{.1in}

\noindent {\bf Theorem 13.} $IFS^-(\forall L(V_{\alpha}))\ \rightarrow\
\aleph_{\alpha}= \aleph_{\alpha}^{L(V_{\omega+\alpha})}\ \wedge\
L(V_{\omega+\alpha+1})\models |V_{\omega+\alpha}|=\aleph_{\alpha}$\\
\proof By induction on $\alpha.$ If $\alpha$ is a limit ordinal then
certainly for $\beta<\alpha$ we have by the induction hypothesis,
$$\aleph_{\beta}=\aleph_{\beta}^{L(V_{\omega+\beta})}$$
which implies
$\aleph_{\alpha}=\aleph_{\alpha}^{L(V_{\omega+\alpha})}.$ By the
induction hypothesis we also have
$|V_{\omega+\alpha}|=\aleph_{\alpha}.$ So in $V$ there is a subset of
$V_{\omega+\alpha}\times V_{\omega+\alpha}$ which is a well ordering
of $V_{\omega+\alpha}$ of order type $\aleph_{\alpha}.$ Since
$P(V_{\omega+\alpha}\times V_{\omega+\alpha})\in
L(V_{\omega+\alpha+1}),$
$$L(V_{\omega+\alpha+1})\models |V_{\omega+\alpha}|=\aleph_{\alpha}$$
If $\alpha=\beta+1$ then since by the induction hypothesis we have
that $\aleph_{\beta}=\aleph_{\beta}^{L(V_{\omega+\beta})}$ and
$L(V_{\omega+\alpha})\models |V_{\omega+\beta}|=\aleph_{\beta},$  all
order types of ordinals less than $\aleph_{\alpha}$ are 
incoded by subsets of $V_{\omega+\beta}\times V_{\omega+\beta}.$ 
Since $P(V_{\omega+\beta}\times V_{\omega+\beta})\in
L(V_{\omega+\alpha}),$ this implies
$\aleph_{\al}=\aleph_{\al}^{L(V_{\omega+\alpha})}.$   Let $P$ be the partial
order of all one to one maps from initial segments of 
$\aleph_{\alpha}$ into $V_{\omega+\alpha}.$ $P$ is $\aleph_{\beta}$
closed and $P$ is a definable element of $L(V_{\omega+\al}).$ 
By $IFS^-(\forall L(V_{\alpha})),$ there exists a $G\subseteq
P$ which is $L(V_{\omega+\alpha})$ generic. $\bigcup G$ is a bijection
from $\aleph_{\alpha}$ onto $V_{\omega+\alpha}.$ Since there is a
subset of $V_{\omega+\alpha}\times V_{\omega+\alpha}$ which is a well
ordering $V_{\omega+\alpha}$ of order type
$\aleph_{\alpha}$ and $P(V_{\omega+\alpha}\times V_{\omega+\alpha})\in
L(V_{\omega+\alpha+1}),$ $V_{\omega+\alpha}$ has size
$\aleph_{\alpha}$ in $L(V_{\omega+\alpha+1}).$ 

\vspace{.1in}

\noindent {\bf Corollary 14.} $IFS^-(\forall L(V_{\alpha}))\vdash GCH$

\vspace{.1in}

\noindent {\bf Theorem 15.} $IFS(\forall L(V_{\alpha}))\ \rightarrow\
\aleph_{\alpha}= \aleph_{\alpha}^{L(V_{\omega+\alpha})}\ \wedge\
L(V_{\omega+\alpha+1})\models |V_{\omega+\alpha}|=\aleph_{\alpha}$\\
\proof Exactly the same as for $IFS^-(\forall L(V_{\alpha})).$

\vspace{.1in}

\noindent {\bf Corollary 16.} $IFS(\forall L(V_{\alpha}))\vdash
IFS^-(\forall L(V_{\alpha})).$

\vspace{.1in}

\noindent {\bf Theorem 17.} 
$IFS^-(\forall L(V_{\al}))$ implies that for every regular cardinal
$\aleph_{\al},$ 
$$L(V_{\omega+\al+1})-L(V_{\omega+\al})\neq\emptyset$$
\proof Suppose not. Let $\aleph_{\al}$ be the least regular cardinal
such that 
$$L(V_{\omega+\al+1})=L(V_{\omega+\al})$$
Note that $\al$ is definable in $L(V_{\omega+\al})$ either as the
least $\be$ such that $V=L(V_{\omega+\be})$ or as the
least $\be$ such that $V=L(V_{\omega+\be-1}).$ Let $P$ be the set of
bijections between subsets of $\aleph_{\al+1}$ of size less than
$\aleph_{\al}$ into subsets of $\aleph_{\al}.$ $P$ is a definable
element of $L(V_{\omega+\al+1})$ and since
$L(V_{\omega+\al+1})=L(V_{\omega+\al}),$ $P$ is a definable element of 
$L(V_{\omega+\al}).$ As $\aleph_{\al}$ is regular, $P$ is
$<\aleph_{\al}$ closed. By $IFS^-(\forall L(V_{\al})),$ there is a
$L(V_{\omega+\al})$ generic subset $G$ of $P$ in $V.$ $\bigcup G$ is a
bijection from $\aleph_{\al+1}$ to $\aleph_{\al}$ a contradiction.

\vspace{.1in}

\noindent {\bf Theorem 18.} $CIFS$ implies $IFS^-(\forall L(V_{\alpha})).$\\
\proof Let $\aleph_{\alpha}$ be regular and $P$ definable in
$L(V_{\omega+\alpha})$ such that $P$ is $\aleph_{\beta}$ closed for
every $\beta<\alpha.$ Since $P$ is definable in $L(V_{\omega+\alpha})$
so is the collection ${\cal D}$ of dense subsets of $P$ in
$L(V_{\omega+\alpha})$ so $|{\cal D}|\leq \aleph_{\alpha}.$ List
${\cal D}$ as $\big\{D_{\zeta}\mid \zeta<\aleph_{\alpha}\big\}.$ Now by
induction on $\zeta<\aleph_{\alpha},$ by the $<\aleph_{\alpha}$
closedness of $P$ we can build a sequence $\langle p_{\zeta}\mid
\zeta<\aleph_{\alpha}\rangle$ such that $p_{\zeta}\in D_{\zeta}.$ Let
$G$ be the filter generated by the $\langle p_{\zeta}\mid
\zeta<\aleph_{\alpha}\rangle.$ Now let $\aleph_{\alpha}$ be singular
and $P$ definable in $L(V_{\omega+\alpha})$ such that $P$ is
$\aleph_{\beta}$ closed for every $\beta<\alpha.$ Since
$\aleph_\alpha$ is singular, $P$ is also $\aleph_{\alpha}$ closed. $P$
is definable in $L(V_{\omega+\alpha+1})$ and so is the set ${\cal D}$
of dense subsets of $P$ in $L(V_{\omega+\alpha+1}),$ so $|{\cal
D}|\leq \aleph_{\alpha+1}.$ As before we can build an
$L(V_{\omega+\alpha+1})$ generic subset of $P.$ 

\vspace{.1in}

\noindent {\bf Theorem 19.} $CIFS\ 
\rightarrow\ $ for each regular cardinal $\aleph_{\alpha},$  every set
definable in $L(V_{\omega+\alpha})$ has size a most $\aleph_{\alpha}$
in $L(V_{\omega+\alpha+1}).$\\
\proof Let $x$ be definable in $L(V_{\omega+\al})$ where
$\aleph_{\al}$ is regular. By theorem 17 $\al$ is definable in
$L(V_{\omega+\al}).$ Let $\ga$ be the least ordinal greater than $\aleph_{\al}$ such that $x\in
L_{\ga}(V_{\omega+\al}).$ Since $\ga$ is definable in
$L(V_{\omega+\al}),$ $\gamma$ has size $\aleph_{\al}.$ Therefore
$(L_{\ga}(V_{\omega+\al}),\in)$ is 
isomorphic to a model $(V_{\omega+\al},E)$ where $E$ is a subset of
$V_{\omega+\al}\times V_{\omega+\al}.$ Since the Mostowski Collapsing
Theorem holds in 
$L(V_{\omega+\alpha+1}),$ $(V_{\alpha+\omega},E)$ is isomorphic
to $(L(V_{\omega+\alpha}),\in)$  in
$L(V_{\omega+\alpha+1})$ and therefore $x$ can have size at most
$|V_{\omega+\al}|=\aleph_{\al}$ in $L(V_{\omega+\alpha+1}).$

\pagebreak

\begin{center}
Picture of the Universe under $ZFC+CIFS$
\end{center}

\vspace{6 in}

\noindent Under $CIFS,$ for every regular cardinal $\aleph_{\al},$ all
definable elements of $L(V_{\omega+\alpha})$ have size at most
$\aleph_{\al}$ in $L(V_{\omega+\alpha+1}),$ forcing
$L(V_{\omega+\al+1})-L(V_{\omega+\al})$ large.

\pagebreak

\begin{center}
{\bf SOME PARTING PHILOSOPHICAL REMARKS}
\end{center}

The conventional view of the history of set theory says that Godel in
1938 proved that the consistency of $ZF$ implies the consistency of
$ZFC$ and of $ZFC+GCH,$ and that Cohen with the invention of forcing
proved that $Con(ZF)$ implies $Con(ZF+\neg AC)$ and $Con(ZFC+\neg
GCH),$ but if $IFS(L)$ is correct, a better way to state the history
would be to say that Godel
discovered $L$ and Cohen discovered that there are many generic
extensions of $L.$

The author believes that not all transitive models of $ZFC$ are
created equal and that set theorists should make more active use of
this fact, while they should place less emphasis on relative
consistency results. Some Formalists may object to the Platonistic
slant of this exposition, but a Formalist can always play the game of
pretending to be a Platonist. Finally, the author thinks it is ironic that although mathematics
and especially mathematical
logic is an art noted for its precise and formalized reasoning, 
it seems that in order to solve problems at the frontiers of logic's
foundations we must tackle questions of an esthetic nature of the kind
addressed in this article.

\pagebreak

\begin{center}
REFERENCES
\end{center}

\begin{enumerate}
\item C. C. Chang and J. Keisler, {\em Model Theory},
North Holland Publishing Co.
\item M. Foreman, {\em Potent Axioms}, Transactions of the A.M.S., vol
294 (1986) pp 1-27.
\item C. Freiling, {Axioms of Symmetry: Throwing Darts at the Real
Line}, this {\sc Journal}, vol. 51 (1988) pp 190-200.
\item {[Jech1] T. Jech, {\em Multiple Forcing}, Cambridge
University Press.}
\item {[Jech2] T. Jech, {\em Set Theory}, Academic Press.}
\item {[Kunen] K. Kunen, {\em Set Theory}, Studies in Logic and the Foundations
of Mathematics, vol 102 (1980), Elsevier Science Publishing Company,
Amsterdam.} 
\item {S. Mac Lane, {\em Is Mathias an Ontologist?}, in Set Theory of
the Continuum, H. Judah, W. Just, and H. Woodin editors, Springer
Verlag (1992) pp 119-122}
\item {[BA I] P. Maddy, {\em Believing the Axioms I}, this {\sc Journal} vol 53
(1988) pp 481-511.}
\item {[BA II] P. Maddy, {\em Believing the Axioms II}, this {\sc Journal} vol 54
(1988) pp 736-764.}
\item {A. R. D. Mathias, {\em What is Mac Lane Missing?}, in Set Theory of
the Continuum, H. Judah, W. Just, and H. Woodin editors, Springer
Verlag (1992) pp 113-118}
\item {R. Penrose, {\em The Emperors New Mind}, Oxford University
Press, Oxford (1989)}
\end{enumerate}

\end{document}